\documentclass{article}
\usepackage{graphicx,amssymb}
\usepackage{graphics}
\usepackage{epsfig}
\usepackage{subfigure}
\usepackage{wrapfig}
\usepackage[T1]{fontenc}
\usepackage[latin9]{inputenc}
\begin{document}


\title{A Generalization of Prefactored Compact Schemes for Advection Equations\\
}

\author{Adrian Sescu  \\
Mississippi State University, Mississippi State, MS 39762, USA \vspace{8mm} \\
\textit{Tel}: 1-662-325-7484, fax: 1-662-325-7730 \\
\textit{Email}: sescu@ae.msstate.edu}

\maketitle

\begin{abstract}

A generalized prefactorization of compact schemes aimed at reducing the stencil and improving the computational efficiency is proposed here in the framework of transport equations. By the prefactorization introduced here, the computational load associated with inverting multi-diagonal matrices is avoided, while the order of accuracy is preserved. The prefactorization can be applied to any centered compact difference scheme with arbitrary order of accuracy (results for compact schemes of up to sixteenth order of accuracy are included in the study). One notable restriction is that the proposed schemes can be applied in a predictor-corrector type marching scheme framework. Two test cases, associated with linear and nonlinear advection equations, respectively, are included to show the preservation of the order of accuracy and the increase of the computational efficiency of the prefactored compact schemes.

\end{abstract}



\section{Introduction}

Compact difference schemes, as opposed to explicit schemes, possess the advantage of attaining higher-order of accuracy with fewer grid points per stencil. They are preferred in applications where high accurate results are desired, such as direct numerical simulations, large eddy simulations, computational aeroacoustics or electromagnetism, to enumerate few, and in some instances they feature accuracy comparable to spectral methods. One of the disadvantages of compact schemes is that an implicit approach is required to determine the derivatives, where a matrix (usually multi-diagonal) has to be inverted.

A comprehensive study of high-order compact schemes approximating both first and second derivatives on a uniform grid was performed by Lele \cite{lele}. A wavenumber based optimization was introduced wherein the dispersion error was reduced significantly (assuming an exact temporal integration). Over the next years, compact schemes have been studied by many research groups, and applied to various engineering problems (see for example, Li et al. \cite{Li}, Adams and Shariff \cite{Adams}, Liu \cite{Liu}, Deng and Maekawa \cite{Deng1}, Fu and Ma \cite{Fu1,Fu2}, Meitz and Fasel \cite{Meitz}, Shen et al. \cite{Shen}, Shah et al. \cite{Shah}). Other examples include Kim and Lee \cite{Kim} who performed an analytic optimization of compact finite difference schemes, Mahesh \cite{Mahesh} who derived a family of compact finite difference schemes for the spatial derivatives in the Navier-Stokes equations based on Hermite interpolations (see also, Chu and Fan \cite{Chu} for a similar prior analysis), or Deng and Zhang \cite{Deng2} who developed compact high-order nonlinear schemes which are equivalent to fifth-order upwind biased explicit schemes in smooth regions. Hixon and Turkel \cite{Hixon1,Hixon2} derived prefactored high-order compact schemes that use three-point stencils and returns up to eighth-order of accuracy. These schemes combine the tridiagonal compact formulation with the optimized split derivative operators of an explicit MacCormack type scheme. The optimization of Hixon's \cite{Hixon1,Hixon2} schemes in terms of reducing the dispersion error was performed by Ashcroft and Zhang \cite{Ashcroft} who used Fourier analysis to select the coefficients of the biased operators such that the dispersion characteristics match those of the original centered compact scheme and their numerical wavenumbers have equal and opposite imaginary components. Sengupta et al. \cite{Sengupta} derived a new compact schemes for parallel computing. Today, compact schemes are widely used in numerical simulations of turbulent flows (e.g., direct numerical simulations), computational aeroacoustics, or computational electromagnetics. In order to increase the speed of such numerical simulations it is desirable to derive more computational efficient compact schemes without affecting the order of accuracy and the wavenumber characteristics.

In this work, we propose a generalized prefactorization of existing compact schemes aimed at reducing the stencil and increasing the computational efficiency. It is based on the type of prefactorization introduced previously by Hixon and Turkel \cite{Hixon1,Hixon2}, but here the order of accuracy can be increased indefinitely, and there are no specific requirements for the original compact schemes to be suitable to prefactorization, other than they must fall in the class of centered scheme. A similar optimization was recently performed by Bose and Sengupta \cite{Bose}. They developed an alternate direction bidiagonal (ADB)
scheme, which showed neutral stability and good dispersion characteristics. The analysis included here can be viewed as a generalization of their work, except the focus is on the order of accuracy rather the dispersion characteristics. One of the restrictions of the schemes derived here is that the proposed prefactored schemes can be combined with a predictor-corrector type time marching scheme only. This allows the determination of the derivatives by sweeping from one boundary to the other, thus avoiding the inversion of matrices which can make the computational algorithm cumbersome and the execution time-consuming. It is shown that the original order of accuracy of the classical compact schemes is preserved, while the computational efficiency can be almost doubled (numerical tests pertaining to fourth and sixth order accurate schemes show over $40\%$ decrease in the processing time, while higher order schemes are expected to be more efficient). 

Section II will discuss the classical compact difference schemes including wavenumber characteristics. In section III, the prefactored compact schemes are introduced and analyzed. In section IV, two test cases are considered to verify the efficiency and to check the preservation of the order of accuracy of the proposed schemes; the two cases correspond to a linear and a nonlinear problem. The last section is reserved for concluding remarks.

\section{Compact Difference Schemes}


Consider the general compact centered approximation for the first derivative:

\begin{equation}\label{e1}
\sum_{k=1}^{N_c}\alpha_k(u_{j+k}' + u_{j-k}') + u_{j}' =
\frac{1}{h}
\sum_{k=1}^{N_e}a_k(u_{j+k} - u_{j-k}) 
 + O(h^n),
\end{equation}
where $1 \le j \le N$ (with $N$ being the number of grid points), the gridfunctions at the nodes are $u_j = u(x_j)$, the values of the derivatives with respect to $x$ are $u_j'$, and $h$ is the spatial step. If $\alpha_k = 0$ then the scheme is termed explicit. Compact schemes (also known as implicit or Pade schemes), by contrast, have $\alpha_k \ne 0$ and require the solution of a matrix equation to determine the derivatives of the grid function. Conventionally, the coefficients $\alpha_k$ and $a_k$ are chosen to give the largest possible exponent, $n$, in the truncation error, for a given stencil width. Table~\ref{t1} gives several weights for various compact difference schemes that are considered in this study: fourth (C4), sixth (C6), eight (C8), tenth (C10), twelfth (C12), fourteenth (C14) and sixteenth (C16) order accurate. Schemes C4 and C6 require a tri-diagonal matrix inversion, C8 and C10 a penta-diagonal matrix inversion, C12 and C14 a seven-diagonal matrix inversion, while C16 a nine-diagonal matrix inversion. The inversion of three- and five-diagonal matrices for determining the weights of compact differencing schemes received increased attention, while seven- and nine-diagonal matrices are less popular due to the computational inefficiency. The prefactored compact scheme of Hixon~\cite{Hixon2} is also included here in the form:

\begin{eqnarray}\label{e2}
a u_{j+1}^{F'} + c u_{j-1}^{F'} + (1-a-c) u_{j}^{F'} =
\frac{1}{h}
\left[
b u_{j+1} - (2b-1) u_j - (1-b) u_{j-1}) 
\right], \nonumber \\
c u_{j+1}^{B'} + a u_{j-1}^{B'} + (1-a-c) u_{j}^{B'} =
\frac{1}{h}
\left[
(1-b) u_{j+1} - (2b-1) u_j - b u_{j-1}) 
\right],
\end{eqnarray}
where $F$ and $B$ stand for 'forward' and 'backward', respectively. For sixth order accuracy, $a=1/2-1/(2\sqrt{5})$, $b=1-1/(30 a)$ and $c=0$.

\begin{table}[htpb]
 \begin{center}
  \caption{Weights of several compact difference schemes}
  \label{t1}
  \begin{tabular}{rrrrrrrrrr} \hline
       Scheme & $\alpha_1$ & $\alpha_2$ & $\alpha_3$ & $\alpha_4$ & $a_1$ & $a_2$ & $a_3$ & $a_4$    \\\hline
       $C4$ &  1/4 &  0 &  0 &  0 &  3/4 & 0 & 0 & 0 \\
       $C6$ &  1/3 &  0 &  0 &  0 &  7/9 & 1/36 & 0 & 0 \\
       $C8$ &  4/9 &  1/36 &  0 &  0 &  20/27 & 25/216 & 0 & 0 \\
       $C10$ &  1/2 &  1/20 &  0 &  0 &  17/24 & 101/600 & 1/600 & 0 \\
       $C12$ &  9/16 &  9/100 &  1/400 &  0 &  21/32 & 231/1000 & 49/4000 & 0 \\
       $C14$ &  3/5 &  3/25 &  1/175 &  0 &  31/50 & 67/250 & 283/12250 & 1/9800 \\
       $C16$ &  16/25 &  4/25 &  16/1225 &  1/4900 &  72/125 & 38/125 & 1784/42875 & 761/686000 \\
\hline
  \end{tabular}
 \end{center}
\end{table}

\begin{figure}[htp]
\begin{center}
      \label{f1a}\includegraphics[width=0.48\textwidth]{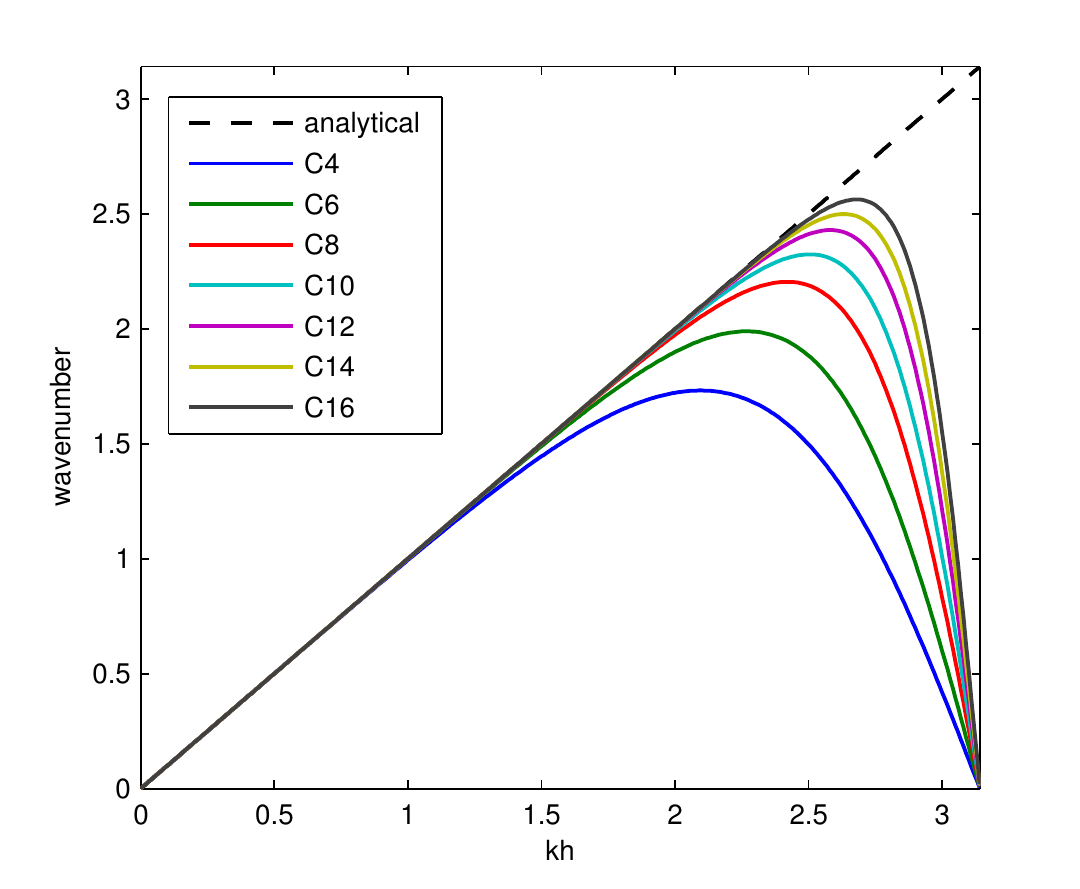}
      \label{f1b}\includegraphics[width=0.48\textwidth]{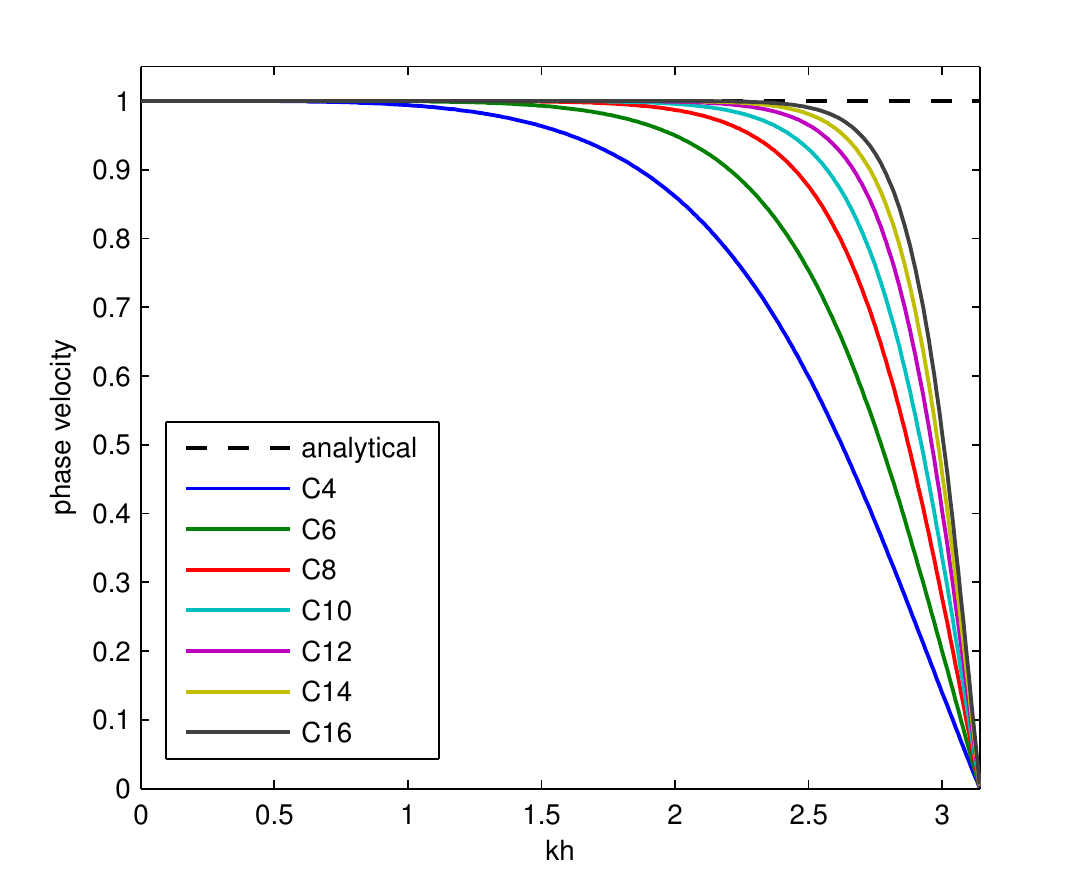}
 \hspace{26mm}(a) \hspace{73mm} (b)
 \end{center}
\begin{center}
      \label{f1b}\includegraphics[width=0.48\textwidth]{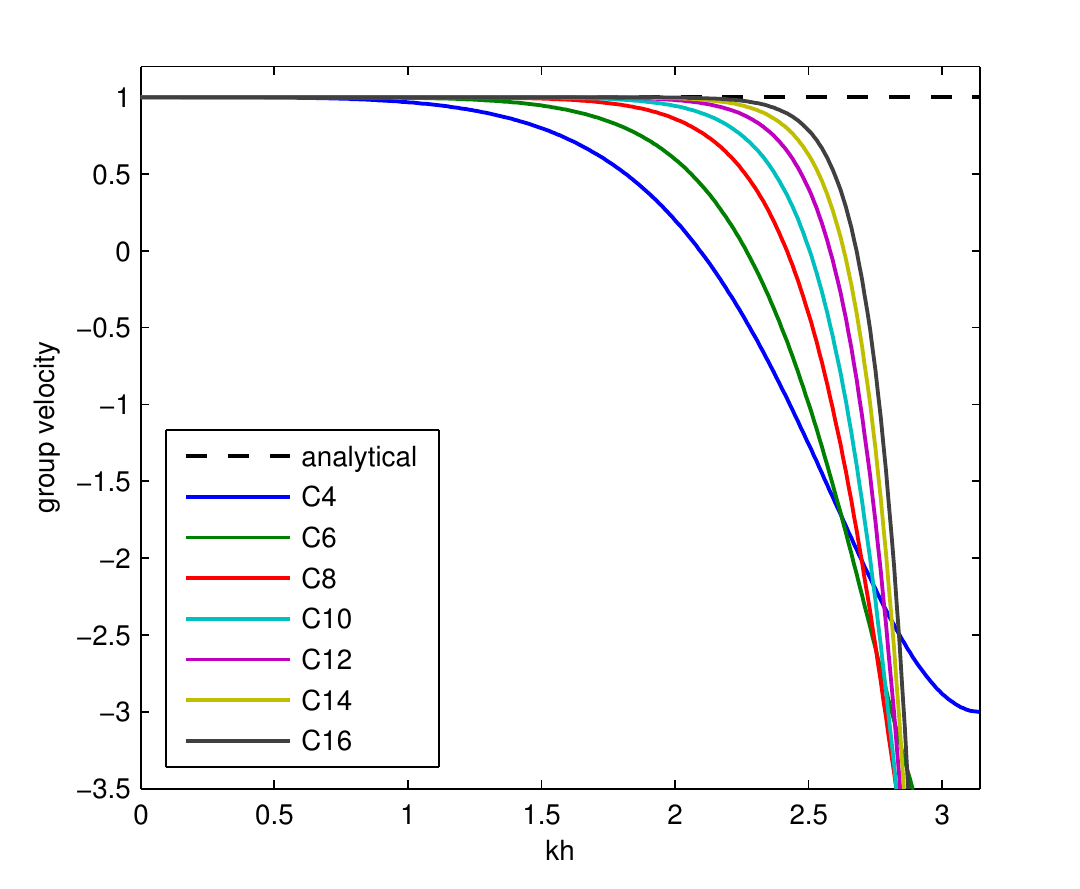}
\hspace{73mm} (c)
 \end{center}
  \caption{a) Numerical wavenumber for various compact schemes compared to the exact wavenumber; b) Numerical phase velocity compared to the exact phase velocity; c) Numerical group velocity compared to the exact group velocity.}
  \label{f1}
\end{figure}

The leading order term in the truncation error of a finite difference scheme depends on the choice of the coefficients and the $(n+1)$st derivative of the function $u$. To obtain the wavenumber characteristics of compact schemes, consider a periodic domain with $N$ uniformly spaced points on $x \in [0,L]$ (with $h=L/N$), and the discrete Fourier transform of $u$ as

\begin{equation}\label{}
u_j =
\sum_{m=-N/2}^{N/2-1} \hat{u}_m e^{ik_m x_j}; \hspace{6mm}
j=1,...,N.
\end{equation}
where the wavenumber is $k_m = 2\pi m/L$. The $m$th component of the discrete Fourier transform of $u'$ denoted $\hat{u}_m'$ is simply $ik_m \hat{u}_m$. Taking the discrete Fourier transform of equation (\ref{e1}) provides the approximate value of $\hat{u}_m'$ in the form:

\begin{equation}\label{}
(\hat{u}_m')_{num} =
iK(k_m h)\hat{u}_m,
\end{equation}
where the numerical wavenumber is given by:

\begin{equation}\label{}
K(z) =
\frac
{\sum_{n=1}^{N_e}2a_n \sin{(nz)}}
{1+\sum_{n=1}^{N_c}2\alpha_n \cos{(nz)}}.
\end{equation}
where $z$ is the wavenumber. Figures \ref{f1} shows the numerical wavenumber, the phase velocity, and the group velocity corresponding to the schemes given in table \ref{t1} (the exact wavenumber, phase, and group velocity are also included for comparison). One can notice that the dispersion error decreases as the order of accuracy is increased.

\section{Generalized profactored schemes}

The classical compact difference scheme (\ref{e1}) is written (conveniently) in the following form:

\begin{equation}\label{ll}
\sum_{k=1}^{N_c}\gamma_k(u_{j+k}' + u_{j-k}') + \left( 1 - 2 \sum_{k=1}^{N_c}\gamma_k  \right) u_{j}' =
\frac{1}{h}
\sum_{k=1}^{N_e}\eta_k(u_{j+k} - u_{j-k}) 
 + O(h^n),
\end{equation}
where the coefficients $\gamma_k$ and $\eta_k$ are functions of the original coefficients $\alpha_k$ and $a_k$ given in table \ref{t1}. The prefactored compact schemes proposed in this work are given in the form

\begin{equation}\label{dd1}
\left( 1 -  \sum_{k=1}^{N_c}\beta_k  \right) u_{j}^{'F} 
+ \sum_{k=1}^{N_c}\beta_k u_{j+k}^{'F} =
\frac{1}{h}
\sum_{k=1}^{N_e} b_k \left( u_{j+k} - u_{j} \right)
+ O(h^n),
\end{equation}
for the forward operator, and

\begin{equation}\label{dd2}
\left( 1 -  \sum_{k=1}^{N_c}\beta_k  \right) u_{j}^{'B} 
+ \sum_{k=1}^{N_c}\beta_k u_{j-k}^{'B} =
\frac{1}{h}
\sum_{k=1}^{N_e} b_k \left( u_{j} - u_{j-k} \right)
+ O(h^n),
\end{equation}
for the backward operator ('forward' and 'backward' correspond to the predictor and corrector steps, respectively). Notice that the two new schemes are of downwind and upwind types, respectively, so they may feature dissipation errors. However, when combining the predictor and corrector operators the dissipation errors equate to zero, while the dispersion error is the same as the one corresponding to the original centered compact scheme. If we consider a spatial discretization on a one-dimensional grid (consisting of $N$ grid points), the schemes (\ref{dd1}) and (\ref{dd2}) can be written in matricial form (Hixon and Turkel \cite{Hixon1,Hixon2}) as

\begin{equation}\label{tt1}
B_{mn}^F U_{n}^{'F}= \frac{1}{h} C_{mn}^F U_{n}
\end{equation}

\begin{equation}\label{tt2}
B_{mn}^B U_{n}^{'B} = \frac{1}{h} C_{mn}^B U_{n}
\end{equation}
where $\{U'\}$ is the vector of derivatives, $\{U\}$ is the vector of grid functions, and $m,n \in \{ 1,2,...,N \}$. Equations (\ref{dd1}) and (\ref{dd2}) imply that

\begin{equation}\label{v1}
B_{mn}^F  = B_{nm}^B
\end{equation}

\begin{equation}\label{v2}
C_{mn}^F  = -C_{nm}^B
\end{equation}
The coefficients in (\ref{dd1}) and (\ref{dd2}) are determined by requiring that the original classical compact scheme is recovered by performing an average between the predictor and corrector operators, formally written as 

\begin{equation}\label{qq}
\left< \cdot \right>^{'} = \frac{1}{2} \left( \left< \cdot \right>^{'F} + \left< \cdot \right>^{'B} \right)
\end{equation}
Multiplying (\ref{qq}) by $B_{mn}^BB_{nm}^B$ and using (\ref{tt1}), (\ref{tt2}), (\ref{v1}) and (\ref{v2}), as well as the relation $B_{mn}^B B_{nm}^B = B_{nm}^B B_{mn}^B$, which is true for matrices of the type considered here, we obtain

\begin{equation}\label{kk}
B_{mn}^BB_{nm}^B\left< \cdot \right>^{'} = \frac{1}{2h} \left( B_{nm}^B C_{mn}^B 
                                                                                              - B_{mn}^B C_{nm}^B  \right)\left< \cdot \right>
\end{equation}

The coefficients $\beta_k$ and $b_k$ can now be determined by matching equation (\ref{kk}) with equation (\ref{ll}) (in the appendix, a Mathematica \cite{math} code used to determine the coefficients for the new 12th order accurate scheme is included). Tables \ref{t2} and \ref{t3} include these coefficients for several schemes of different orders of accuracy, fourth (PC4), sixth (PC6), eigth (PC8), tenth (PC10), twelfth (PC12), fourteenth (PC14) and sixteenth (PC16), that correspond to the classical compact schemes given in table \ref{t1}.

\begin{table}[htpb]
 \begin{center}
  \caption{Weights of several prefactored compact difference schemes (left-hand-side).}
  \label{t2}
  \begin{tabular}{rrrrrrrrrr} \hline
       Scheme & $\beta_1$ & $\beta_2$ & $\beta_3$ & $\beta_4$     \\\hline
       $PC4$ &    0.211324870586 &  0 &  0 &  0  \\
       $PC6$ &    0.276393202250 &  0 &  0 &  0  \\
       $PC8$ &    0.353614989057 &  0.022913166676 &  0 &  0  \\
       $PC10$ &  0.390891054882 &  0.041982762456 &  0 &  0  \\
       $PC12$ &  0.424261339307 &  0.076528671307 &  0.002177424900 &  0  \\
       $PC14$ &  0.440844836186 &  0.103628733678 &  0.005175974177 &  0  \\
       $PC16$ &  0.450833811211 &  0.139274137394 &  0.012291382216 &  0.000195518547  \\
\hline
  \end{tabular}
 \end{center}
\end{table}

\begin{table}[htpb]
 \begin{center}
  \caption{Weights of several prefactored compact difference schemes (right-hand-side).}
  \label{t3}
  \begin{tabular}{rrrrrrrrrr} \hline
       Scheme & $b_1$ & $b_2$ & $b_3$ & $b_4$    \\\hline
       $PC4$ &    1.000000000000 & 0 & 0 & 0 \\
       $PC6$ &    0.907868932583 & 0.046065533708 & 0 & 0 \\
       $PC8$ &    0.679849926548 & 0.160075036725 & 0 & 0 \\
       $PC10$ &    0.544199349631 & 0.223702048938 & 0.002798850830 & 0 \\
       $PC12$ &    0.377436479527 & 0.283739040458 & 0.018361813185 & 0 \\
       $PC14$ &    0.270368050633 & 0.312589794656 & 0.034570978390 & 0.000184856220 \\
       $PC16$ &    0.157403729700 & 0.326389389050 & 0.060796869771 & 0.001856720721 \\
\hline
  \end{tabular}
 \end{center}
\end{table}

The stencil count has been reduced as follows: from 3- to 2-point stencil for PC4 (this is similar to the scheme proposed by Hixon \cite{Hixon2}), from 5- to 3-point stencil for PC6 and PC8, from 7- to 4-point stencil for PC10 and PC12, and from 9- to 5-point stencil for PC14 and PC16. The advantage of these prefactored schemes is that there is no need to invert matrices because the derivatives can be obtained explicitly by sweeping from one boundary to the other (assuming the grid functions and the derivatives are available at the boundaries). This simplifies the computational cost significantly, without affecting the performance of the schemes since by averaging the predictor and corrector operators the original classical compact schemes are obtained. Thus, the wavenumber characteristics manifested through zero-dissipation and low dispersion of the original compact scheme (\ref{e1}) is retained.

\section{Test cases}

\subsection{Preliminaries}

We consider the initial-value problem in $\textbf{R}\times [0,\infty)$:

\begin{equation}\label{ww}
\frac {\partial{u}}{\partial{t}}
+c(u) \frac {\partial{u}}{\partial{x}} = 0,
\end{equation}

\begin{equation}\label{}
u(x,0)=u_0(x),
\end{equation}
and appropriate boundary conditions, where $u(x,t)$ is a scalar function, $c$ is the convective velocity which may depend on $u$, and $u_0(x)$ is a given function of space. Let $\Omega=\{ x,l_1<x<l_2 \}$ be a finite domain in the real set $\textbf{R}$ with $l_1$ and $l_2$ chosen such that there exist a real non-negative number $h=(l_1-l_2)/N$ called spatial step ($N$ is an integer representing the number of grid points).

Because the accuracy of the time marching scheme is not the focus of this study, we use a second order MacCormack~\cite{MacCormack} scheme which is a two-step predictor-corrector time advancement scheme, and a second order TVD Runge-Kutta scheme~\cite{liu2}. The first time marching schemes is applied in the framework of prefactored compact schemes, while the second scheme is applied in the framework of classical compact schemes. To increase the accuracy of the time marching scheme, a very small time step is set, the emphasis being on the accuracy of the spatial discretization. For the advection equation in one dimension, the MacCormack scheme can be written as

\begin{eqnarray}\label{}
u_{j}^F   &=& u_{j}^{n} - \sigma \Delta_{x}^{F} u_{j}^{n}  \nonumber \\
u_{j}^B  &=& u_{j}^F    - \sigma \Delta_{x}^{B} u_{j}^F \\
u_{j}^{n+1} &=& \frac{1}{2} (u_{j}^F + u_{j}^B), \nonumber
\end{eqnarray}
where $\sigma = c(u) k/h $, $k$ is the time step, $\Delta_{x}^{F}$ is the 'forward', and $\Delta_{x}^{B}$ is the 'backward' difference operators. The second order TVD Runge-Kutta method~\cite{liu2} is:

\begin{eqnarray}\label{37}
u_{j}^{(0)} &=& u_{j}^n                 \nonumber \\
u_{j}^{(1)} &=& u_{j}^{(0)}+\Delta t L(u_{j}^{(0)}) \\
u_{j}^{n+1} &=& \frac{1}{2}u_{j}^{(0)}+\frac{1}{2}u_{j}^{(1)}+\frac{1}{2}\Delta t L(u_{j}^{(2)}) \nonumber
\end{eqnarray}
where $L(u_{j}) = -c(u_{j}) (\partial{u}/\partial{x})_{j}$. These two time marching schemes are essentially the same, except the first is applied in the framework of predictor-corrector type marching procedure, while the second is applied with 2 stages per time step.

\subsection{Linear advection equation}

\begin{figure}[!htpb]
 \begin{center}
      \includegraphics[width=15cm]{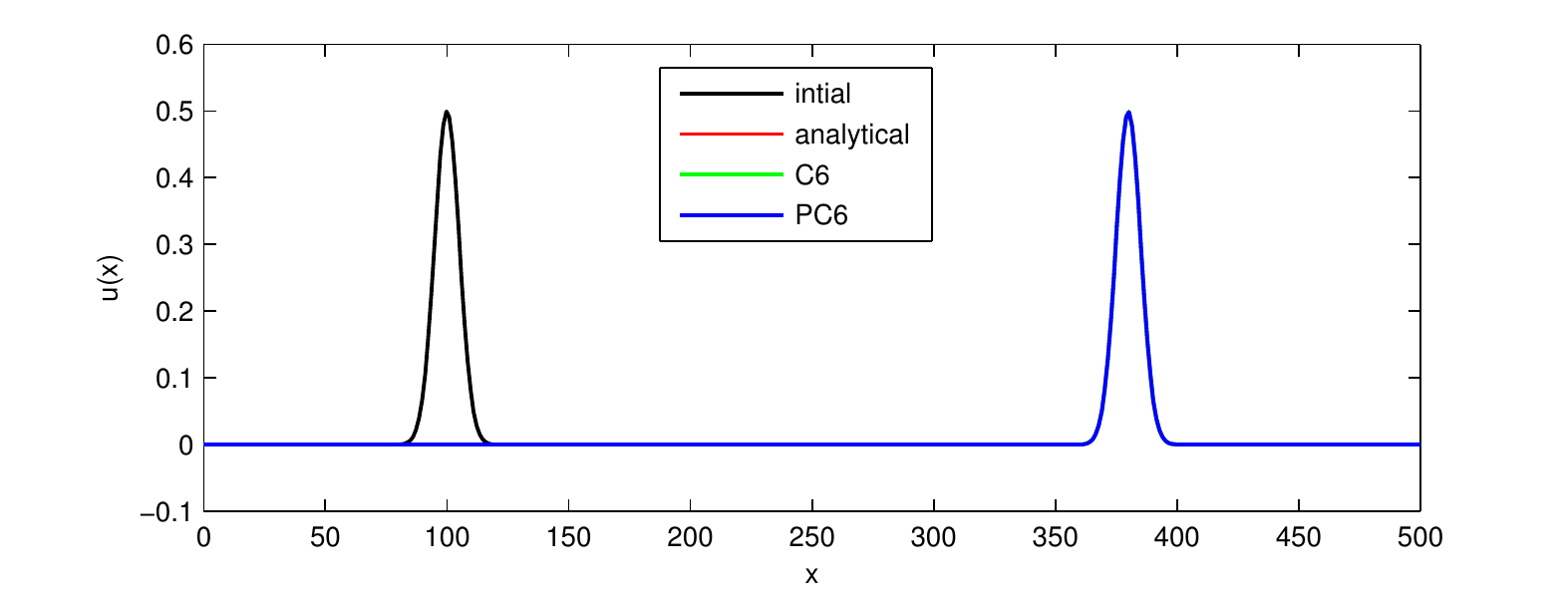}
 \end{center}
  \caption{Plot of the linear solution at $t=280$ for the sixth order accurate schemes.}
  \label{f2}
\end{figure}

For the linear advection equation, the convective velocity $c$ is a positive constant (equal to $1$ here), which renders the transport of the initial condition $u_0(x)$ in the positive direction. The initial condition is given by

\begin{eqnarray}\label{}
u_0(x) = u(x,0) = \frac{1}{2} exp\left[ -(\ln 2) \frac{x^2}{9} \right]
\end{eqnarray}
(Hardin et al. \cite{Hardin}). The domain boundaries are $l_1=-20$ and $l_2=450$, and the final time is $t_f=200$. The equation (\ref{ww}) is solved numerically using all prefactored compact difference schemes considered in the previous section. Figure \ref{f2} shows the numerical solution using the sixth order accurate prefactored and classical compact schemes; the numerical solutions are compared to the analytical solution (the initial solution is also included). One can notice that there is no difference among the solutions. This is stressed in table \ref{t8} which lists the $L_2$-errors for the fourth, sixth, eight and tenth order accurate schemes (both prefactored and classical); the differences between the errors are very small, which is expected since the schemes have similar behavior (the differences in the errors come from the time marching algorithms).

\begin{table}[htpb]
 \begin{center}
  \caption{\footnotesize $L_2$ error comparison between prefactored and conventional compact schemes.}
  \label{t8}
  \begin{tabular}{rrrrrr} \hline
    Order of accuracy & $L_2$ error $PC4$ & $L_2$ error $C4$ \\
\hline
    4th      &  4.4090580981E-003     &  4.4101995527E-003   \\
    6th      &  6.0846046295E-004     &  6.0895634631E-004  \\
    8th      &  7.2394864471E-005     &  7.2375087716E-005  \\
    10th    &  1.5989324725E-005     &  1.5470197664E-005  \\
\hline
  \end{tabular}
 \end{center}
\end{table}

\begin{figure}[!htpb]
 \begin{center}
    \mbox{
      \subfigure {\label{f2_a}\includegraphics[width=6.6cm]{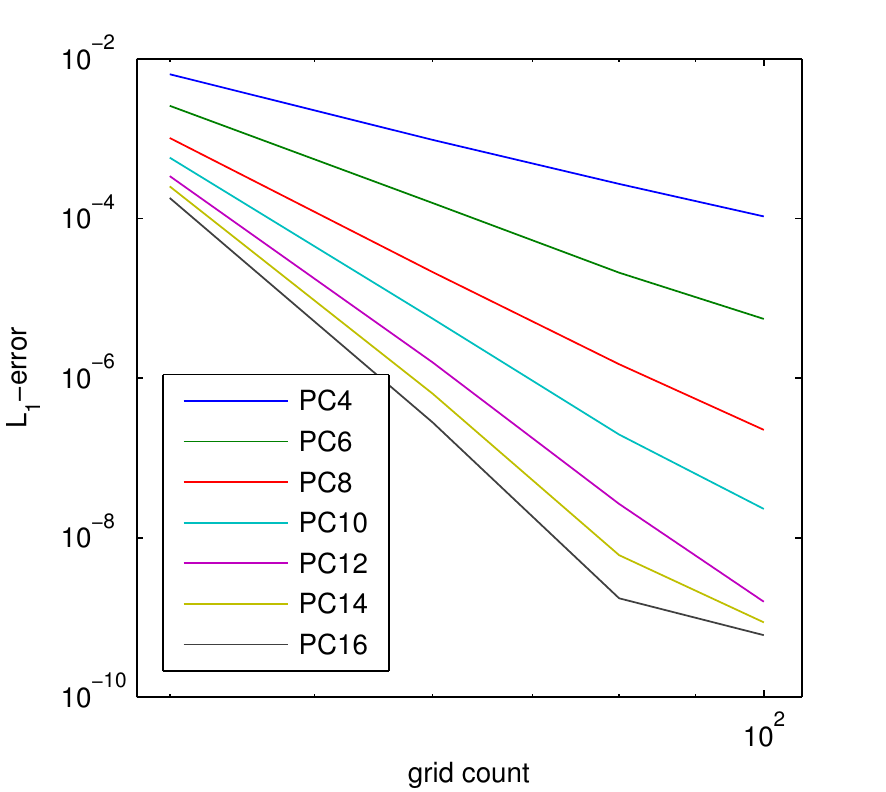}}
      \subfigure {\label{f2_b}\includegraphics[width=6.6cm]{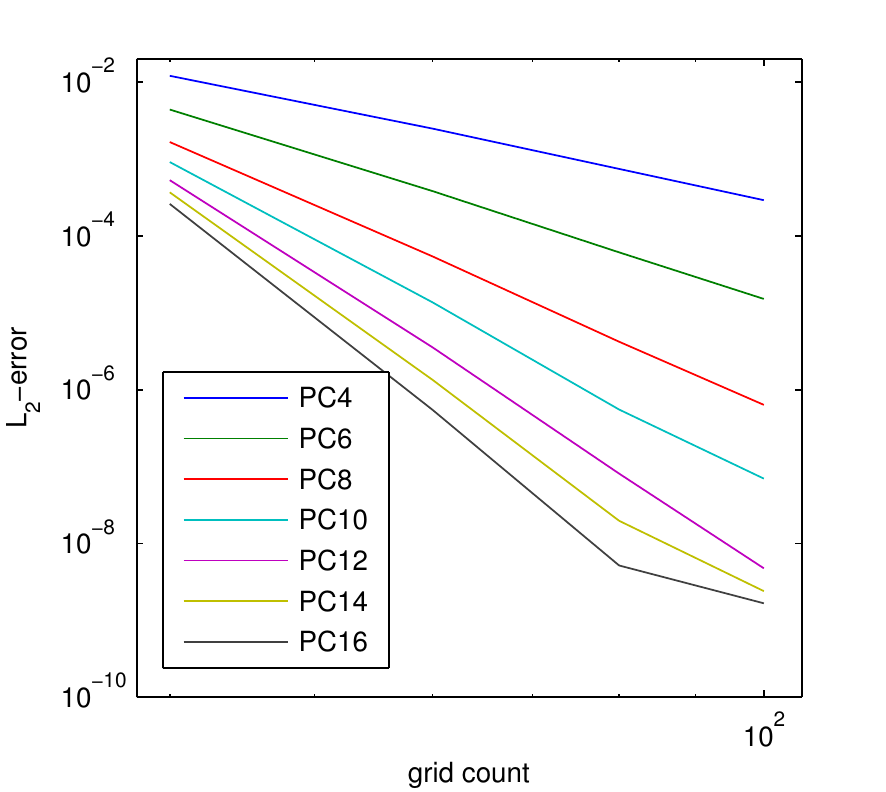}} 
      }\\
a) \hspace{70mm} b)
 \end{center}
 \begin{center}
    \mbox{
      \subfigure {\label{f2_c}\includegraphics[width=6.6cm]{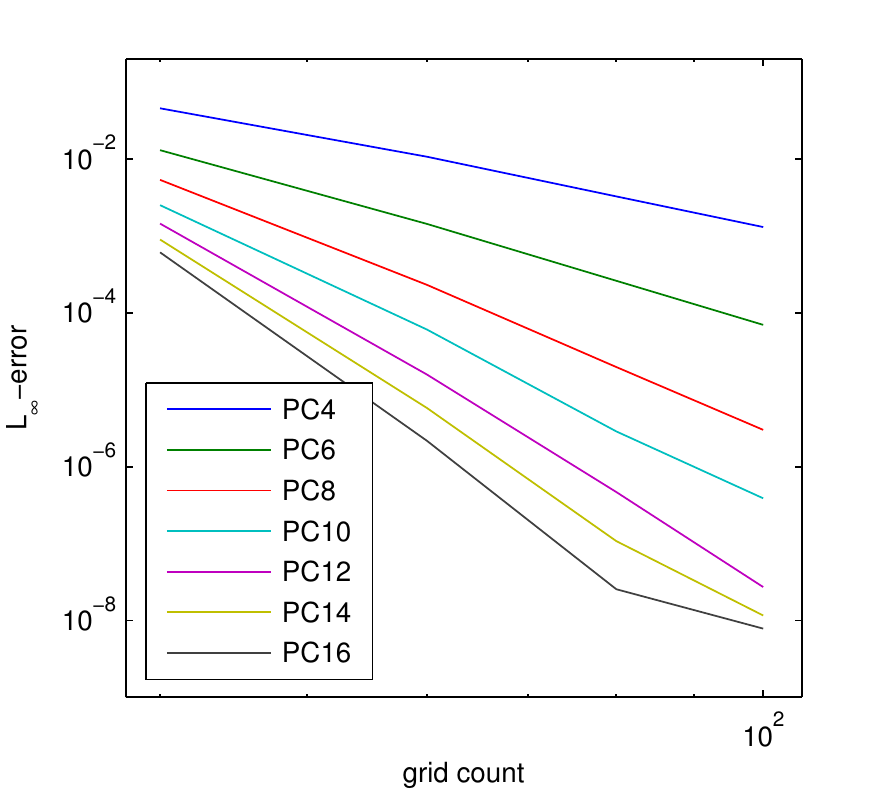}} 
      }\\
      c)
 \end{center}
  \caption{$L_1$, $L_2$ and $L_{\infty}$ errors plotted against grid point count.}
  \label{f3}
\end{figure}

From Taylor series expansions corresponding to two different steps, $h_1$ and $h_2$, the order of accuracy of the compact scheme can be estimated as

\begin{eqnarray}\label{qqq}
p = \frac{\ln(\epsilon_1/\epsilon_2)}{\ln(h_1/h_2)}
\end{eqnarray}
where $\epsilon_1$ and $\epsilon_2$ are errors associated with the spatial steps $h_1$ and $h_2$, respectively. The smallest time step in the MacCormack time marching scheme was $k = 1E-6$ corresponding the the sixteenth order accurate scheme. Different grid point counts are considered to study the behavior of the numerical errors. Figure \ref{f3} shows $L_1$, $L_2$ and $L_{\infty}$ errors (calculated by taking a summation over all points in the grid) plotted against grid point count, corresponding to all prefactored schemes, while table \ref{t4} lists these errors and the evaluated order of accuracy (calculated via equation (\ref{qqq})). From the plots in figure \ref{f3}, a linear decrease of all errors with respect to the grid count can be noticed, except slight deviations for the last grid count in the behavior of PC14 and PC16. At this accuracy level the errors approach the machine precision, but in addition the time marching may no longer provide the required accuracy (a further decrease of the time step was not pursued). The evaluated orders of accuracy included in table \ref{t4} are close to the theoretical ones (indicated by the suffixes, e.g. 'PC4' has 4th order of accuracy) suggesting that the prefactored schemes are indeed able to maintain the desired order.

The computational efficiency increase is tested for PC4, PC6, PC8 and PC10, by comparing the computation time to the case when corresponding classical compact schemes of fourth (C4), sixth (C6), eight (C8) and tenth (C10) order of accuracy, respectively, are employed. C4 and C6 schemes necessitate the inversion of a three-diagonal matrix which is performed here via the Thomas algorithm (\cite{Thomas}), while for the C8 and C10 schemes a pentadiagonal matrix inversion is required; this is done here using a LU-factorization, where the matrix elements are calculated only once and stored to be re-used. For classical compact schemes, a second order TVD Runge-Kutta scheme - which is equivalent to the second-order MacCormack schemes in terms of the number of operations - is used to march the solution in time. The number of grid points is significantly large - in the order of 10,000 - such that the spatial discretization consumes most of the computational time. Table \ref{t7} shows the results in terms of the computational time. By using PC4 as opposed to C4, a $41.7\%$ decrease in the computational time is achieved; a $40.2\%$ decrease in the computational time is achieved when employing PC6 scheme versus C6; a $33.1\%$ decrease is achieved when applying PC6 versus C6; and a $32.3\%$ decrease is achieved when applying PC10 versus C10. 

The percentages shown in table \ref{t7} can be explained by a comparison in terms of the number of operations that are necessary when the two types of schemes are applied. It was found that fewer operations are necessary in the case of PCn schemes in all cases (for example, 4 additions and 5 multiplications per step are needed for 8th order accurate PCn schemes, as opposed to 7 additions and 8 multiplications per stage that are needed for 8th order accurate Cn schemes)

\begin{table}[htpb]
 \begin{center}
  \caption{$L_1$, $L_2$ and $L_{\infty}$ errors for different grid point counts - linear case.}
  \label{t4}
  PC4 scheme
  \begin{tabular}{rrrrrrrrrr} \hline
       Grid & $L_1$-error & $L_2$-error  & $L_{\infty}$-error      \\\hline
       $40$ &    6.4309202E-003   &   1.2147772E-002   &   4.5647113E-002 \\
       $60$ &    9.7376776E-004   &   2.5019901E-003   &   1.0738430E-002 \\
       $80$ &    2.7148411E-004   &   7.4323530E-004   &   3.2702186E-003 \\
       $100$ &  1.0678123E-004   &   2.9341239E-004   &   1.3099043E-003 \\
\hline
  \end{tabular}
  \begin{tabular}{rrrrrr}
      estimated order of accuracy, p = 3.9058  \\  \hline
  \end{tabular}
 \end{center}
 \begin{center}
  PC6 scheme
  \begin{tabular}{rrrrrrrrrr} \hline \hspace{2.3mm} 
       $40$ &    2.5853741E-003   &   4.4073207E-003   &   1.3058909E-002 \\
       $60$ &    1.5631156E-004   &   3.8415139E-004   &   1.4335771E-003 \\
       $80$ &    2.1072380E-005   &   6.1152769E-005   &   2.6060378E-004 \\
       $100$ &  5.5345770E-006   &   1.5190237E-005   &   6.9914708E-005 \\
\hline
  \end{tabular}
  \begin{tabular}{rrrrrr}
      estimated order of accuracy, p = 5.8700  \\  \hline
  \end{tabular}
 \end{center}
 \begin{center}
  PC8 scheme
  \begin{tabular}{rrrrrrrrrr} \hline \hspace{2.3mm} 
       $40$ &    1.0269107E-003   &   1.6601837E-003   &   5.3893200E-003 \\
       $60$ &    2.1429783E-005   &   5.3936887E-005   &   2.3061682E-004 \\
       $80$ &    1.4878953E-006   &   4.2014630E-006   &   1.9661071E-005 \\
       $100$ &  2.2419349E-007   &   6.3432171E-007   &   3.0175207E-006 \\
\hline
  \end{tabular}
  \begin{tabular}{rrrrrr}
      estimated order of accuracy, p = 8.3155  \\  \hline
  \end{tabular}
 \end{center}
 \begin{center}
  PC10 scheme
  \begin{tabular}{rrrrrrrrrr} \hline \hspace{2.3mm} 
       $40$ &    5.7962247E-004   &   9.1838126E-004   &   2.5118507E-003 \\
       $60$ &    5.5528859E-006   &   1.3631097E-005   &   6.0412057E-005 \\
       $80$ &    1.9686060E-007   &   5.5050250E-007   &   2.8709444E-006 \\
       $100$ &  2.2838208E-008   &   6.9681216E-008   &   3.8666113E-007 \\
\hline
  \end{tabular}
  \begin{tabular}{rrrrrr}
      estimated order of accuracy, p = 10.1140  \\  \hline
  \end{tabular}
 \end{center}
 \begin{center}
  PC12 scheme
  \begin{tabular}{rrrrrrrrrr} \hline \hspace{2.3mm} 
       $40$ &    3.3956164E-004   &   5.3364060E-004   &   1.4421609E-003 \\
       $60$ &    1.5771299E-006   &   3.5518382E-006   &   1.5715408E-005 \\
       $80$ &    2.6383987E-008   &   8.0943016E-008   &   4.6658106E-007 \\
       $100$ &  8.7747738E-009   &   2.4115746E-008   &   1.1512857E-007 \\
\hline
  \end{tabular}
  \begin{tabular}{rrrrrr}
      estimated order of accuracy, p = 11.9781  \\  \hline
  \end{tabular}
 \end{center}
 \begin{center}
  PC14 scheme
  \begin{tabular}{rrrrrrrrrr} \hline \hspace{2.3mm} 
       $40$ &    2.5287271E-004   &   3.6913391E-004   &   9.0116478E-004 \\
       $60$ &    6.3975306E-007   &   1.3428221E-006   &   5.8091922E-006 \\
       $80$ &    6.0074747E-009   &   1.9783115E-008   &   1.0758384E-007 \\
       $100$ &  2.8940807E-009   &   7.8821369E-009   &   3.7617685E-008 \\
\hline
  \end{tabular}
  \begin{tabular}{rrrrrr}
      estimated order of accuracy, p = 13.4793  \\  \hline
  \end{tabular}
 \end{center}
 \begin{center}
  PC16 scheme
  \begin{tabular}{rrrrrrrrrr} \hline \hspace{2.3mm} 
       $40$ &    1.8127733E-004   &   2.6222753E-004   &   6.1433898E-004 \\
       $60$ &    2.7869625E-007   &   5.4486475E-007   &   2.1557661E-006 \\
       $80$ &    5.6610923E-009   &   1.5792774E-008   &   7.7338168E-008 \\
       $100$ &  2.8717976E-009   &   7.8132634E-009   &   3.7871620E-008 \\
\hline
  \end{tabular}
  \begin{tabular}{rrrrrr}
      estimated order of accuracy, p = 15.2329  \\  \hline
  \end{tabular}
 \end{center}
\end{table}

\begin{table}[htpb]
 \begin{center}
  \caption{\footnotesize Computational time decrease by using prefactored schemes.}
  \label{t7}
  \begin{tabular}{rrrrrr} \hline
    Schemes & comp. time decrease ($\%$) \\
\hline
    $PC4$ vs. $C4$     &  41.7\%   \\
    $PC6$ vs. $C6$     &  40.2\%  \\
    $PC8$ vs. $C8$     &  33.1\%  \\
    $PC10$ vs. $C10$ &  32.3\%  \\
\hline
  \end{tabular}
 \end{center}
\end{table}

\subsection{Nonlinear advection equation}

\begin{figure}[!htpb]
 \begin{center}
    \mbox{
      \subfigure {\label{f2_a}\includegraphics[width=7.5cm]{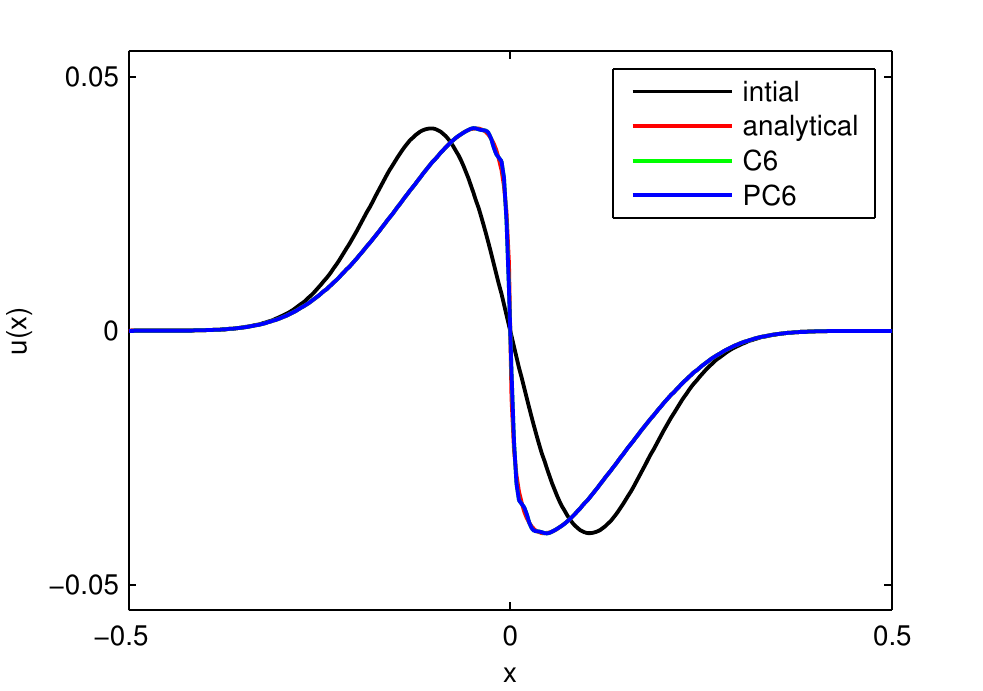}}
      \subfigure {\label{f2_b}\includegraphics[width=7.5cm]{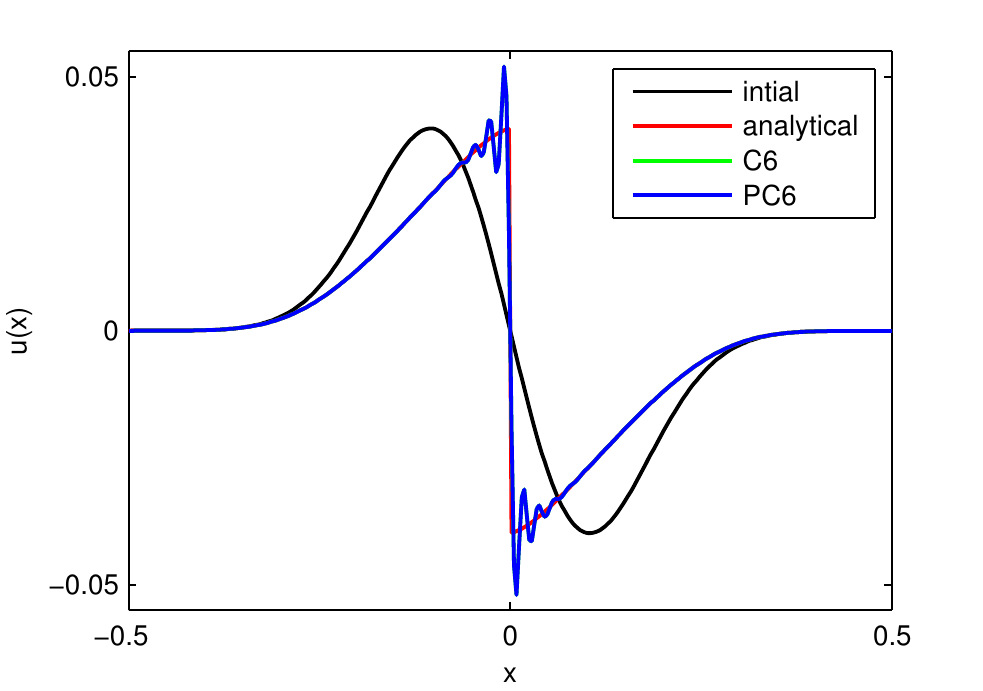}} 
      }\\
a) \hspace{70mm} b)
 \end{center}
  \caption{Plots of the nonlinear solution for the sixth order accurate schemes: a) $t=1.6$; b) $t=2.2$.}
  \label{f4}
\end{figure}

In the context of the inviscid Burgers equation, the convective velocity is $c(u) = u$, which renders equation (\ref{ww}) to be nonlinear. The domain boundaries are $l_1 = -1/2$ and $l_2 = 1/2$, while the initial condition is given as

\begin{eqnarray}\label{}
u_0(x) = u(x,0) = 0.1exp\left( -\frac{x^2}{0.16^2} \right) \sin(2\pi x)
\end{eqnarray}
where the Gaussian function was introduced to drive the initial waveform exponentially to zero, at both boundaries (in this ways, errors from boundary conditions are minimized). The final time is $t_f=1.5$. As in the linear case, the equation (\ref{ww}) is solved numerically using all prefactored compact difference schemes for different grid point counts. Figure \ref{f4} shows the numerical solution at two different time instances, using the sixth order accurate prefactored and classical compact schemes; a comparison to the analytical solution is also included along with the initial condition. A discontinuity is forming in $x=0$, which is stronger in figure \ref{f4}a, corresponding to $t=2.2$; here both compact schemes develop spurious oscillations around the discontinuity, but the important conclusion is that both prefactored and conventional schemes behave similarly, as expected. Table \ref{t8} lists the $L_2$-errors for the fourth, sixth, eight and tenth order accurate schemes (both prefactored and classical schemes); the differences between the errors are very small as in the linear case.

\begin{table}[htpb]
 \begin{center}
  \caption{\footnotesize $L_2$ error comparison between prefactored and conventional compact schemes.}
  \label{t10}
  \begin{tabular}{rrrrrr} \hline
    Order of accuracy & $L_2$ error $PC4$ & $L_2$ error $C4$ \\
\hline
    4th      &  1.9978793015E-007     &  2.0023910853E-007   \\
    6th      &  7.9774807530E-008     &  7.8986651180E-008  \\
    8th      &  3.1242097731E-008     &  2.8410542257E-008  \\
    10th    &  1.4655545221E-008     &  5.7531058004E-009  \\
\hline
  \end{tabular}
 \end{center}
\end{table}

\begin{figure}[!htpb]
 \begin{center}
    \mbox{
      \subfigure {\label{f2_a}\includegraphics[width=6.6cm]{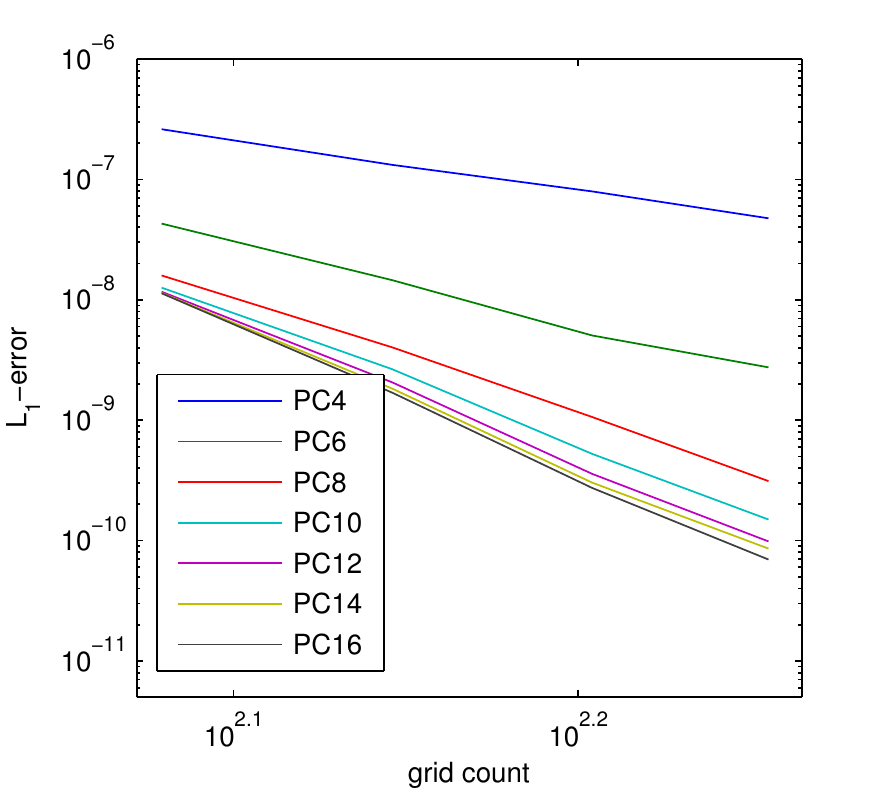}}
      \subfigure {\label{f2_b}\includegraphics[width=6.6cm]{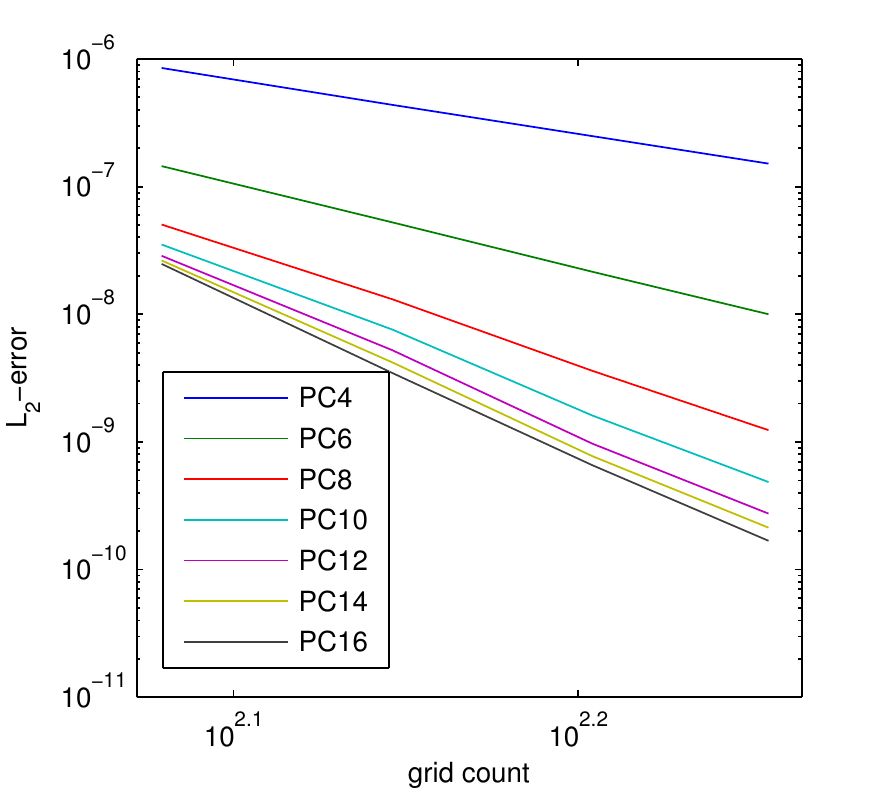}} 
      }\\
a) \hspace{70mm} b)
 \end{center}
 \begin{center}
    \mbox{
      \subfigure {\label{f2_c}\includegraphics[width=6.6cm]{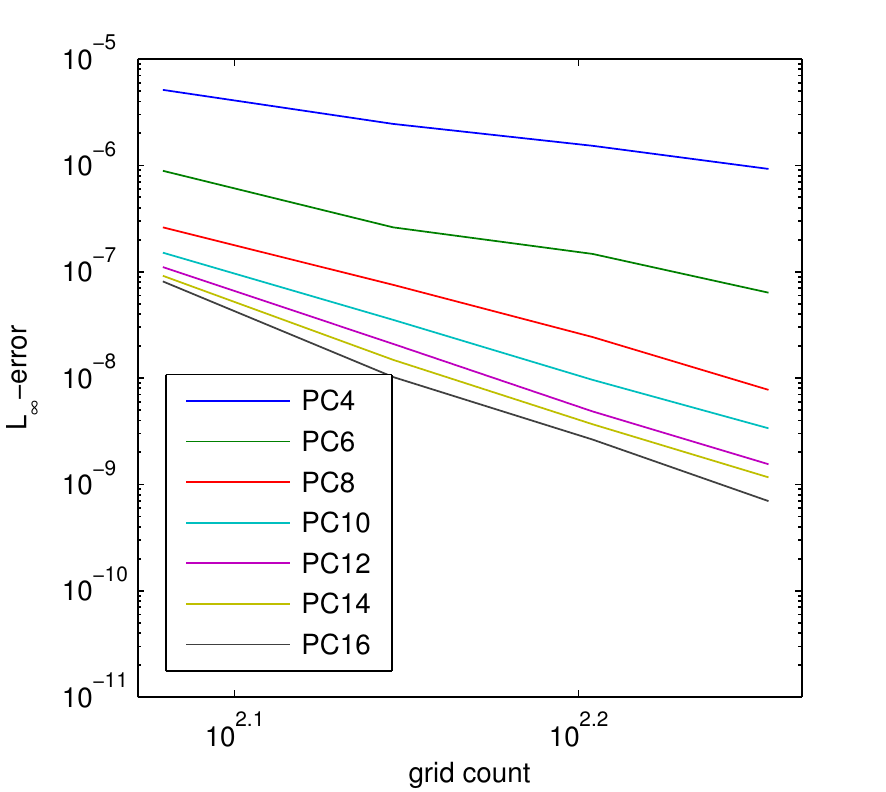}} 
      }\\
      c)
 \end{center}
  \caption{$L_1$, $L_2$ and $L_{\infty}$ errors plotted against grid point count.}
  \label{f5}
\end{figure}

Figure \ref{f5} shows $L_1$, $L_2$ and $L_{\infty}$ errors plotted against grid point count, corresponding to all prefactored schemes. Up to tenth order of accuracy, the trends are the same as in the linear case, but for order of accuracy greater than twelfth there are some discrepancies that may be the result of the accuracy of calculating the exact solution. Table \ref{t5} lists $L_1$, $L_2$ and $L_{\infty}$ errors and the evaluated order of accuracy (calculated via equation (\ref{qqq})). The orders of accuracy included in table \ref{t5} are close to the theoretical ones for PC4, PC6, PC8 and PC10, but there are deviations for PC12, PC14 and PC16 (due to the same reasons that are mentioned above). However, the trend of increasing the evaluated order of accuracy as the theoretical one is increased is still right.

In terms of the computational efficiency gain, by using PC4 as opposed to C4 a $40.6\%$ decrease in the computational time is achieved; a $39.8\%$ decrease in the computational time is achieved when employing PC6 scheme versus C6; a $31.6\%$ decrease is achieved when applying PC8 scheme versus C8; and a $30.4\%$ decrease is achieved when applying PC8 scheme versus C8 (see table \ref{t9}).

\begin{table}[htpb]
 \begin{center}
  \caption{$L_1$, $L_2$ and $L_{\infty}$ errors for different grid point counts - nonlinear case.}
  \label{t5}
  PC4 scheme
  \begin{tabular}{rrrrrrrrrr} \hline
       Grid  & $L_1$-error & $L_2$-error  & $L_{\infty}$-error      \\\hline
       $120$ &    2.6057545E-007   &   8.5267800E-007   &   5.1480039E-006 \\
       $140$ &    1.3189979E-007   &   4.3851923E-007   &   2.4635989E-006 \\
       $160$ &    7.9331220E-008   &   2.4853181E-007   &   1.5294059E-006 \\
       $180$ &    4.7629938E-008   &   1.5153784E-007   &   9.2540099E-007 \\
\hline
  \end{tabular}
  \begin{tabular}{rrrrrr}
      estimated order of accuracy, p = 4.2243  \\  \hline
  \end{tabular}
 \end{center}
 \begin{center}
  PC6 scheme
  \begin{tabular}{rrrrrrrrrr} \hline \hspace{2.3mm} 
       $120$ &    4.2923461E-008   &   1.4477527E-007   &   8.8853969E-007 \\
       $140$ &    1.4589376E-008   &   5.2417825E-008   &   2.6191106E-007 \\
       $160$ &    5.0682594E-009   &   2.1581203E-008   &   1.4697422E-007 \\
       $180$ &    2.7495252E-009   &   1.0043626E-008   &   6.3695264E-008 \\
\hline
  \end{tabular}
  \begin{tabular}{rrrrrr}
      estimated order of accuracy, p = 6.3415  \\  \hline
  \end{tabular}
 \end{center}
 \begin{center}
  PC8 scheme
  \begin{tabular}{rrrrrrrrrr} \hline \hspace{2.3mm} 
       $120$ &    1.5900337E-008   &   5.0506253E-008   &   2.6150391E-007 \\
       $140$ &    4.0357681E-009   &   1.3133096E-008   &   7.5480377E-008 \\
       $160$ &    1.0600384E-009   &   3.6076941E-009   &   2.4301185E-008 \\
       $180$ &    3.1257899E-010   &   1.2368444E-009   &   7.7387924E-009 \\
\hline
  \end{tabular}
  \begin{tabular}{rrrrrr}
      estimated order of accuracy, p = 8.4347  \\  \hline
  \end{tabular}
 \end{center}
 \begin{center}
  PC10 scheme
  \begin{tabular}{rrrrrrrrrr} \hline \hspace{2.3mm} 
       $120$ &    1.2630357E-008   &   3.5129113E-008   &   1.5106339E-007 \\
       $140$ &    2.6550133E-009   &   7.5831916E-009   &   3.5363402E-008 \\
       $160$ &    5.2494587E-010   &   1.6113340E-009   &   9.6153273E-009 \\
       $180$ &    1.4983937E-010   &   4.8509852E-010   &   3.3684929E-009 \\
\hline
  \end{tabular}
  \begin{tabular}{rrrrrr}
      estimated order of accuracy, p = 9.8317  \\  \hline
  \end{tabular}
 \end{center}
 \begin{center}
  PC12 scheme
  \begin{tabular}{rrrrrrrrrr} \hline \hspace{2.3mm} 
       $120$ &    1.1642690E-008   &   2.8799170E-008   &   1.1068219E-007 \\
       $140$ &    2.0583881E-009   &   5.2232619E-009   &   2.1042827E-008 \\
       $160$ &    3.5757610E-010   &   9.6967151E-010   &   4.8609743E-009 \\
       $180$ &    9.8160882E-011   &   2.7427907E-010   &   1.5507093E-009 \\
\hline
  \end{tabular}
  \begin{tabular}{rrrrrr}
      estimated order of accuracy, p = 11.1075  \\  \hline
  \end{tabular}
 \end{center}
 \begin{center}
  PC14 scheme
  \begin{tabular}{rrrrrrrrrr} \hline \hspace{2.3mm} 
       $120$ &    1.1386239E-008   &   2.6357760E-008   &   9.1586511E-008 \\
       $140$ &    1.8347290E-009   &   4.2047321E-009   &   1.4882790E-008 \\
       $160$ &    3.0278882E-010   &   7.7245435E-010   &   3.6903868E-009 \\
       $180$ &    8.5577391E-011   &   2.1319123E-010   &   1.1708572E-009 \\
\hline
  \end{tabular}
  \begin{tabular}{rrrrrr}
      estimated order of accuracy, p = 11.4742  \\  \hline
  \end{tabular}
 \end{center}
 \begin{center}
  PC16 scheme
  \begin{tabular}{rrrrrrrrrr} \hline \hspace{2.3mm} 
       $120$ &    1.1315127E-008   &   2.4784162E-008   &   8.0909695E-008 \\
       $140$ &    1.6970188E-009   &   3.4771499E-009   &   1.0244846E-008 \\
       $160$ &    2.7276821E-010   &   6.5903831E-010   &   2.6510809E-009 \\
       $180$ &    6.9435657E-011   &   1.6853644E-010   &   6.9679461E-010 \\
\hline
  \end{tabular}
  \begin{tabular}{rrrrrr}
      estimated order of accuracy, p = 12.1178  \\  \hline
  \end{tabular}
 \end{center}
\end{table}

\begin{table}[htpb]
 \begin{center}
  \caption{\footnotesize Computational time decrease by using prefactored schemes.}
  \label{t9}
  \begin{tabular}{rrrrrr} \hline
    Schemes & comp. time decrease ($\%$) \\
\hline
    $PC4$ vs. $C4$     &  40.6\%   \\
    $PC6$ vs. $C6$     &  39.8\%  \\
    $PC8$ vs. $C8$     &  31.6\%  \\
    $PC10$ vs. $C10$ &  30.4\%  \\
\hline
  \end{tabular}
 \end{center}
\end{table}

\section{Conclusions}

In this work, a prefactorization of classical compact schemes was proposed, targeting the increase of computational efficiency and the reduction of the stencil. The new prefactored compact schemes can be applied in the context of predictor-corrector type time marching schemes, where sweeping from one boundary to the other is possible. The wavenumber characteristics of the original schemes are preserved since by averaging the predictor and corrector operators the dissipation error vanishes, while the dispersion error matches exactly the dispersion error of the original scheme. Here, up to sixteenth order accurate compact schemes are analyzed, but theoretically any order of accuracy can be considered and tested, including existing optimized compact schemes. As an example, the coefficients of the prefactored spectral-like compact scheme of Lele \cite{lele} were calculated (given in table \ref{t6}) and tested, and the results showed similar behavior as the original scheme.

\begin{table}[htpb]
 \begin{center}
  \caption{Weights of the prefactored spectral-like compact scheme of Lele \cite{lele}.}
  \label{t6}
  \begin{tabular}{rrrrrrrrrr} \hline
        $\beta_1$ & $\beta_2$     \\\hline
         0.4482545282296 &  0.0817278256497  \\
\hline
         $b_1$ & $b_2$ & $b_3$    \\\hline
           0.3069790178973 & 0.3294144889364 & 0.0113973418854 \\
\hline
  \end{tabular}
 \end{center}
\end{table}

$L_1$, $L_2$ and $L_{\infty}$ errors for different grid point counts, both for linear and nonlinear cases, showed that the order of accuracy of the new schemes is preserved (with some deviations in the nonlinear case attributed to other external sources of errors). As expected, the computational efficiency increased by using the prefactored schemes as opposed to corresponding classical compact schemes. This was demonstrated for fourth and sixth order accurate schemes (the extrapolation to higher order scheme seems obvious).

As mentioned above, one of the disadvantages of the proposed prefactored schemes is the restriction to predictor-corrector type time marching schemes. Another disadvantage may be the difficulty in applying the schemes in parallel solvers and at the boundaries of the domain, although similar issues may also be encountered when employing classical high order compact schemes.

\section*{Appendix}

In this appendix, a Mathematica code used to determine the weights of the optimized 12th order accurate scheme is included.

Next matrix A is from Taylor series

A = $\{$$\{ 1, 1, 1, -2, -2, -2 \}$,

$\{ 1, 2^2, 3^2, -2*3! /2!, -2*3!/2!*2^2, -2*3!/2!*3^2\}$,

$\{ 1, 2^4, 3^4, -2*5! /4!, -2*5!/4!*2^4, -2*5!/4!*3^4\}$,

$\{ 1, 2^6, 3^6, -2*7! /6!, -2*7!/6!*2^6, -2*7!/6!*3^6\}$,

$\{ 1, 2^8, 3^8, -2*9! /8!, -2*9!/8!*2^8, -2*9!/8!*3^8\}$,

$\{ 1, 2^{10}, 3^{10}, -2*11!/10!, -2*11!/10!*2^{10}, -2*11!/10!*3^{10}\}$$\}$;

\vspace{4mm}

R = $\{$1,0,0,0,0,0$\}$;

\vspace{4mm}

XF = Inverse[A].R;

\vspace{4mm}

XF[[1]]=XF[[1]]/2; XF[[2]]=XF[[2]]/4; XF[[3]]=XF[[3]]/6;

\vspace{4mm}

Simplify[XF]

\vspace{4mm}

Solve[$\{$a3/(1-2*a1-2*a2-2*a3)=XF[[6]],

a2/(1-2*a1-2*a2-2*a3)=XF[[5]],

a1/(1-2*a1-2*a2-2*a3)=XF[[4]]$\}$,$\{$a1,a2,a3$\}$]

\vspace{4mm}

b1=XF[[1]]*(1-2*a1-2*a2-2*a3)

b2=XF[[2]]*(1-2*a1-2*a2-2*a3)

b3=XF[[3]]*(1-2*a1-2*a2-2*a3)

\vspace{4mm}

NSolve[$\{$ap3*(1-ap1-ap2-ap3)=s3,

ap2*(1-ap1-ap2-ap3)+ap1*ap3=s2,

ap1*(1-ap1-ap2-ap3)+ap1*ap2+ap2*ap3=s1$\}$,

$\{$ap1,ap2,ap3$\}$,WorkingPrecision$->$25]

\vspace{4mm}

par = (1-ap1-ap2-ap3);

\vspace{4mm}

Next 1/2 is from averaging the backward and forward operators;

\vspace{4mm}

NSolve[$\{$1/2(-ae3*par-ap3*(ae1+ae2+ae3))=-b3,

1/2(-ae2*par-ae3*ap1+ap3*ae1-ap2*(ae1+ae2+ae3))=-b2,

1/2(-ae1*par-ae2*ap1-ae3*ap2+af3*ae2+ap2*ae1-ap1*(ae1+ae2+ae3))=-b1$\}$,

$\{$ae1,ae2,ae3$\}$,WorkingPrecision$->$20]



\begin{thebibliography}{9}

 \bibitem{Adams}
 Adams, N.A. and Shariff, K. (1996) A high-resolution hybrid compact-ENO scheme for shock-turbulence interaction problems, {\it J. Comput. Phys.}, Vol. 127, pp. 27.
 
  \bibitem{Ashcroft}
 Ashcroft, G. and Zhang, X. (2003), Optimized prefactored compact schemes, {\it J. Comput. Phys.}, Vol. 190, pp. 459-477.
 
  \bibitem{Bose}
 R. Bose, T.K. Sengupta (2014) Analysis and Design of a New Dispersion Relation Preserving Alternate Direction Bidiagonal Compact Scheme, Journal of Scientific Computing, Vol. 61, pp. 1-28
  
  \bibitem{Bruger}   
 Bruger, A., Gustafsson, B., Lotstedt, P. and Nilsson, J. (2005) High-order accurate solution of the incompressible Navier-Stokes equations, {\it J. Comput. Phys.}, Vol. 203, pp. 49.
  
   \bibitem{Chu}    
  Chu, P.C. and Fan, C. (1998), A three-point combined compact difference scheme, {\it J. Comput. Phys.}, Vol. 140, pp. 370-399.
  
  \bibitem{Deng1}  
 Deng, X., Maekawa, H. (1997) Compact high-order accurate nonlinear schemes, {\it J. Comput. Phys.}, Vol. 130, pp. 77-91.
 
   \bibitem{Deng2} 
 Deng, X. and Zhang, H. (2000), Developing high-order weighted compact nonlinear schemes, {\it J. Comput. Phys.}, Vol. 165, pp. 22-44.
 
  \bibitem{Eswaran}   
 Eswaran, De A.K. (2006), Analysis of a new high resolution upwind compact scheme, {\it J. Comput. Phys.}, Vol. 218, pp. 398-416.

 \bibitem{Fu1}
Fu, D.X. and Ma, Y.W. (1997) A high-order accurate finite difference scheme for complex flow fields, {\it J. Comput. Phys.}, Vol. 134, pp. 1-15.
  
  \bibitem{Fu2}    
 Fu, D., Ma, Y. (2001) Analysis of Super Compact Finite Difference Method and Application to Simulation of Vortex-Shock Interaction, {\it Int. J. Numer. Meth. Fluids}, Vol. 36, pp. 773-805.
 
  \bibitem{Hardin}
 Hardin, J.C., Ristorcelli, J.R. and Tam, C.K.W. (1995), ICASE/LaRC Workshop on benchmark problems in Computational Aeroacoustics, NASA Conference Publication 3300.

 \bibitem{Hixon2}
Hixon, R. (2000) Prefactored small-stencil compact schemes, {\it J. Comput. Phys.}, Vol. 165, pp. 522-541.

 \bibitem{Hixon1}
Hixon, R. and Turkel, E. (2000), Compact implicit MacCormack-type schemes with high accuracy, {\it Journal of Computational Physics}, Vol. 158, pp. 51-70.

 \bibitem{Kim}
Kim, J.W. and Lee, D.J. (1996), Optimized compact finite difference schemes with maximum resolution, {\it AIAA Journal}, Vol. 34, pp. 887-893.

 \bibitem{lele}
 Lele, S.K. (1992) Compact finite difference schemes with spectral-like resolution, {\it Journal of Computational Physics}, Vol. 103, pp. 16-42, 1992.
 
  \bibitem{Li} 
 Li, M., Tang, T. and Fornberg, B. (1995) A compact fourth-order finite difference scheme for the steady incompressible NavierÐStokes equations, {\it Int. J. Numer. Meth. Fluids}, Vol. 20, pp. 1137-1151.
 
  \bibitem{Liu}
 Liu, W. E, J.-G. (1996) Essentially compact schemes for unsteady viscous incompressible flows, {\it J. Comput. Phys.}, Vol. 126, pp. 122-138.
 
 \bibitem{liu2}
 Liu. X., Osher, S. and Chan, T.: {Weighted essentially non-oscillatory schemes}, Journal of Computational Physics 115, 200-212 (1994).

 \bibitem{MacCormack}
MacCormack, Robert W. (1969) The Effect of Viscosity in Hypervelocity Impact Cratering. {\it AIAA Paper AIAA-69-354}.

 \bibitem{Mahesh}
Mahesh, K. (1998), A family of high order finite difference schemes with good spectral resolution, {\it J. Comput. Phys.}, Vol. 145, pp. 332-358.
  
 \bibitem{Meitz}   
 Meitz, H.L. and Fasel, H.F. (2000) A compact difference scheme for the NavierÐStokes equations in velocity-vorticity formulation, {\it J. Comput. Phys.}, Vol. 157, pp. 371.


\bibitem{Sescu2} 
Sescu, A. and Hixon, R. (2013) Numerical anisotropy study of a class of compact schemes, {\it Journal of Scientific Computing}, DOI 10.1007/s10915-014-9826-0.

\bibitem{Sengupta} 
Sengupta, T.K., Dipankara, A. and Kameswara Rao, A. (2007), A new compact scheme for parallel computing using domain decomposition, {\it J. Comput. Phys.}, Vol. 220, pp. 654-677.

\bibitem{Shah} 
Shah, A., Yuan, L. and Khan, A. (2010) Upwind compact finite difference scheme for time-accurate solution of the incompressible NavierÐStokes equations, {\it Applied Mathematics and Computation}, Vol. 215, pp. 3201-3213.

 \bibitem{Shen}
 Shen, Y.Q., Yang, G.W. and Gao, Z. (2006) High-resolution finite compact differences for hyperbolic conservation laws, {\it J. Comput. Phys.}, Vol. 216, pp. 114-137.
  
  \bibitem{Shukla}  
 Shukla, R.K.,  Tatineni, M. and Zhong, X. (2007) Very high-order compact finite difference schemes on non-uniform grids for incompressible NavierÐStokes equations, {\it J. Comput. Phys.}, Vol. 224, pp. 1064-1094.


 \bibitem{Thomas}
Thomas, L.H. (1949), Elliptic Problems in Linear Differential Equations over a Network, Watson Sci. Comput. Lab Report, Columbia University, New York.

 \bibitem{math}
Wolfram Research, Inc. (2014) Mathematica, Version 10.0, Wolfram Research, Inc., Champaign, Illinois.



\end{thebibliography}
\end{document}